\documentclass[reqno,12pt]{amsart}
\usepackage{graphicx,epsfig}
\usepackage{amssymb,amsthm,amsmath,amsfonts,amstext,mathrsfs}
\usepackage{latexsym,url}
\usepackage[T1]{fontenc}
\usepackage{txfonts}
\usepackage[all]{xy}
\newtheorem{thm}{Theorem}
\newtheorem{cor}{Corollary}
\newtheorem{lem}{Lemma}
\newtheorem{prop}{Proposition}

\input xy
\xyoption{all}

\def\N{\hbox{\font\dubl=msbm10 scaled 1200 {\dubl N}}}
\def\Q{\hbox{\font\dubl=msbm10 scaled 1200 {\dubl Q}}}
\def\d{\,{\rm{d}}}

\title[Moments of the question mark function]
{The moments of Minkowski question mark function: the dyadic period function}
\author[Giedrius Alkauskas]{Giedrius Alkauskas}

\begin{document}

\begin{abstract}\rm The Minkowski question mark function $?(x)$ arises as a real distribution of rationals in the Farey
tree. We examine the generating function of moments of $?(x)$. It appears that
the generating function is a direct dyadic analogue of period functions for
Maass wave forms and it is defined in the cut plane
$\mathbb{C}\setminus(0,\infty)$. The exponential generating function satisfies
the integral equation with kernel being the Bessel function. The solution of
this integral equation leads to the definition of dyadic eigenfunctions, arising from a certain Hilbert-Schmidt operator. Finally, we
describe $p-$adic distribution of rationals in the Stern-Brocot tree.
Surprisingly, the Eisenstein series $G_{1}(z)$ does manifest in both real and
$p-$adic cases.
\end{abstract}
\maketitle
\begin{center}
\end{center}
\setcounter{tocdepth}{1}
    \tableofcontents
\section{Introduction}
\begin{figure}[h]
\begin{center}
\includegraphics[width=340pt,height=470pt,angle=-90]{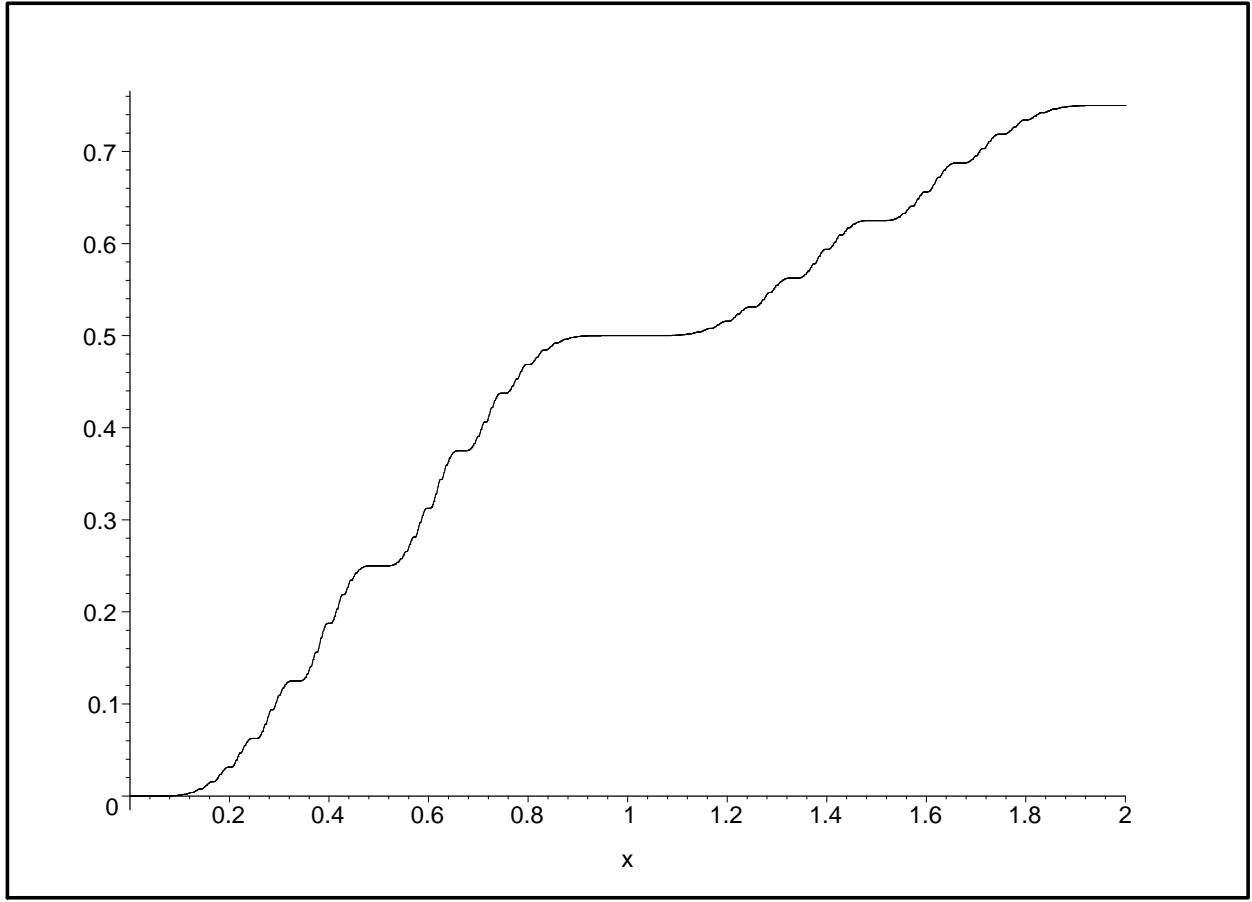}
\caption{The Minkowski question mark function $F(x)$, $x\in[0,2]$}
\end{center}
\end{figure}
\rm This paper\footnote{The current version: October 2008. This is an essential
revision of the first version (September 2006-May 2007)} is the first in the
series of four papers (others being \cite{alkauskas1}, \cite{alkauskas2} and
\cite{alkauskas3}) which are devoted to the study of moments and integral
transforms of
the Minkowski question mark function.\\

The function $?(x)$ (``the question mark function'') was introduced by
Minkowski in 1904 \cite{minkowski} as an example of continuous and monotone
function $?:[0,1]\rightarrow[0,1]$, which maps rationals to dyadic rationals,
and quadratic irrationals to non-dyadic rationals. It is though more convenient
to work with the function $F(x):=?\big{(}\frac{x}{x+1}\big{)}$,
$x\in[0,\infty)\cup\{\infty\}$. Thus, for non-negative real $x$ it is defined
by the expression
\begin{eqnarray}
F([a_{0},a_{1},a_{2},a_{3},...])=1-2^{-a_{0}}+2^{-(a_{0}+a_{1})}-2^{-(a_{0}+a_{1}+a_{2})}+...,
\label{min}
\end{eqnarray}
where $x=[a_{0},a_{1},a_{2},a_{3},...]$ stands for the representation of $x$ by
a (regular) continued fraction \cite{khinchin}. Hence, according to our
convention, $?(x)=2F(x)$ for $x\in[0,1]$. For rational $x$ the series
terminates at the last nonzero element $a_{n}$ of the continued fraction.
 \subsection{Short literature overview}
The Minkowski question mark function was investigated by many authors. In this
subsection we give an overview of available literature.\\

Denjoy \cite{denjoy1} gave an explicit expression for $F(x)$ in terms of
continued fraction expansion; that is, formula (\ref{min}). He also showed that
$?(x)$ is purely singular: the derivative, in terms of the Lebesgue measure,
vanishes almost everywhere. Salem \cite{salem} proved that $?(x)$ satisfies the
Lipschitz condition of order $\frac{\log 2}{2\log \gamma}$, where
$\gamma=\frac{1+\sqrt{5}}{2}$, and this is in fact the best possible exponent
for the Lipschitz condition. The Fourier-Stieltjes coefficients of $?(x)$,
defined as $\int_{0}^{1}e^{2\pi i nx}\d ?(x)$, where also investigated in the
same paper (these coefficients also appeared in \cite{bonanno_isola}; see also
\cite{reese}). The author, as an application of Wiener's theorem about Fourier
series, gives average results on these coefficients without giving an answer to
yet unsolved problem whether these coefficients vanish, as $n\rightarrow\infty$
(it is worth noting that \cite{alkauskas1} analogous Fourier coefficients are
introduced and examined). Kinney \cite{kinney} proved that the Hausdorff
dimension of growth points of $?(x)$ (denote this set by $\mathcal{A}$) is
equal to $\alpha=\frac{1}{2}\big{(}\int_{0}^{1}\log_{2}(1+x)\d ?(x)\big{)}^{-1}$. Numerical estimates for this constant were obtained in
\cite{lagarias2} and \cite{tichy_uitz}; based on the three term functional
equation, we are able to calculate the Kinney's constant to a high precision
(in the appendix of \cite{alkauskas3} we calculate $35$ exact digits; note that
some digits presented in \cite{paradis_viader_bibiloni1} are incorrect). Also,
if $x_{0}\in \mathcal{A}$, $?(x)$ at a point $x_{0}$ satisfies the Lipschitz
condition with exponent $\alpha$. The function $?(x)$ is mentioned in
\cite{conway} in connection with a ``box'' function. In \cite{lagarias_tresser}
Lagarias and Tresser introduce the so called $\mathbb{Q}-$tree: an extension of
the Farey tree which contains all (positive and negative) rationals. Tichy and
Uitz \cite{tichy_uitz} extended Kinney's approach (mainly, calculation of a
Hausdorff dimension) to a parametrized class of singular functions related to
$?(x)$. Bower \cite{bower} considers the solution of the equation $?(x)=x$,
different from $x=0,\frac{1}{2}$ or $1$. There are two of them (symmetric with
respect to $x=\frac{1}{2}$), the first one is given by $x=0.42037233_{+}$
\cite{finch1}. Apparently, no closed form formula exists for it. In
\cite{dilcher_stolarsky} Dilcher and Stolarsky introduced what they call Stern
polynomials. The construction is analogous to similar constructions given in
\cite{esposti_isola_knauf} and \cite{alkauskas3}. Nevertheless, in
\cite{dilcher_stolarsky} all polynomials have coefficients $0$ and $1$, and
their structure is compatible with regular continued fraction algorithm. In
\cite{dushistova_moshchevitin} Dushistova and Moshchevitin find conditions in
order $?'(x)=0$ and $?'(x)=\infty$ to hold (for certain fixed positive real
$x$) in terms of $\limsup_{t\rightarrow\infty}\frac{a_{0}+a_{1}+...+a_{t}}{t}$
and $\liminf_{t\rightarrow\infty}\frac{a_{0}+a_{1}+...+a_{t}}{t}$ respectively,
where $x=[a_{0},a_{1},a_{2},...]$ is represented by a continued fraction. The
nature of singularity of $?(x)$ was clarified by Viader, Parad\'{i}s and
Bibiloni \cite{paradis_viader_bibiloni1}. In particular, the existence of the
derivative $?'(x)$ in $\mathbb{R}$ for fixed $x$ forces it to vanish. Some
other properties of $?(x)$ are demonstrated in \cite{paradis_viader_bibiloni2}.
In \cite{kessebohmer_stratmann2} Kesseb\"{o}hmer and Stratmann studied various
fractal geometric aspects of the Minkowski question mark function $F(x)$. They
showed that the unit interval can be written as the union of three sets:
$\Lambda_{0}:= \{x : F'(x) = 0\}$, $\Lambda_{\infty}:=\{x : F'(x) =\infty\}$,
and $\Lambda_{\sim}:=\{x:F'(x) \text{ does not exist and }F'(x)\neq\infty\}$.
Their main result is that the Hausdorff dimensions of these sets are related in
the following way:
$\dim_{H}(\nu_{F})<\dim_{H}(\Lambda_{\sim})=\dim_{H}(\Lambda_{\infty})=\dim_{H}\big{(}\mathscr{L}(h_{\tt
top})\big{)}<\dim_{H}(\Lambda_{0})=1$. Here $\mathscr{L}(h_{\tt top})$ refers
to the level set of the Stern-Brocot multifractal decomposition at the
topological entropy $h_{\tt top}=\log 2$ of the Farey map $Q$, and
$\dim_{H}(\nu_{F})$ denotes the Hausdorff dimension of the measure of maximal
entropy of the dynamical system associated with $Q$. The notions and technique
were developed earlier by authors in \cite{kessebohmer_stratmann1}.  The paper
\cite{lagarias} deals with the interrelations among the additive continued
fraction algorithm, the Farey tree, the Farey shift and the Minkowski question
mark function. The motivation for the work \cite{panti} is a fact that the
function $?(x)$ can be characterized as the unique homeomorphism of the real
unit interval that conjugates the Farey map with the tent map. In \cite{panti}
Panti constructs an $n$-dimensional analogue of the Minkowski function as the
only homeomorphism of an $n$-simplex that conjugates the piecewise-fractional
map associated to the M\"{o}nkemeyer continued fraction algorithm with an
appropriate tent map. In \cite{bonanno_isola} Bonanno and Isola introduce a
class of $1$-dimensional maps which can be used to generate the binary trees in
different ways, and study their ergodic properties. This leads to studying some
random processes (Markov chains and martingales) arising in a natural way in
this context. In the course of the paper the authors also introduce a function
$\rho(x)=?\big{(}\frac{x}{x+1}\big{)}$, which is, of course, exactly $F(x)$.
Okamoto and Wunsch \cite{okamoto_wunsch} construct yet another generalization
of $?(x)$, though their main concern is to introduce a new family of purely
singular functions. Meanwhile, the paper by Grabner, Kirschenhofer and Tichy
\cite{grabner_tichy}, out of all papers in the bibliography list, is the
closest in spirit to the current article.
 In order to derive precise error bounds for the so called Garcia entropy of a certain measure, the
  authors consider the moments of the continuous and singular function
$
F_{2}([a_{1},a_{2},...])=\sum_{n=1}^{\infty}(-1)^{n-1}3^{-(a_{1}+...+a_{n}-1)}(q_{n}+q_{n-1})$,
where $q_{\star}$ stand for a corresponding denominator of the convergent to
$[a_{1},a_{2},...]$. Lamberger \cite{lamberger} showed that $F(x)$ and
$F_{2}(x)$ are the first two members of a family (indexed by natural numbers)
of mutually singular measures, derived from the subtractive Euclidean
algorithm. The latter two papers are very interesting and promising, and the
author of this article does intend to generalize the results about $F(x)$ to
the whole family $F_{j}(x)$, $j\in\mathbb{N}$.
\subsection{Stern-Brocot tree }Recently, Calkin and Wilf \cite{calkin_wilf} (re-)defined a binary tree which
is generated by the iteration
$$
{a\over b}\quad\mapsto\quad {a\over a+b}\ ,\quad {a+b\over b},
$$
starting from the root ${1\over 1}$ (this tree is a permutation of the
well-known Stern-Brocot tree). Elementary considerations show that this tree
contains any positive rational number once and only once, each represented in
lowest terms \cite{calkin_wilf}. First four iterations lead to
\[
\xymatrix @R=.5pc @C=.5pc { & & & & & & & {1\over 1} & & & & & & & \\
& & & {1\over 2} \ar@{-}[urrrr] & & & & & & & & {2\over 1} \ar@{-}[ullll] & & & \\
& {1\over 3} \ar@{-}[urr] & & & & {3\over 2}\ar@{-}[ull] & & & & {2\over 3}\ar@{-}[urr] & & & & {3\over 1} \ar@{-}[ull] & \\
{1\over 4} \ar@{-}[ur] & & {4\over 3} \ar@{-}[ul] & & {3\over 5}
\ar@{-}[ur] & & {5\over 2} \ar@{-}[ul] & & {2\over 5} \ar@{-}[ur]
& & {5\over 3} \ar@{-}[ul] & & {3\over 4} \ar@{-}[ur] & & {4\over
1} \ar@{-}[ul] }
\]
Thus, the $n$th generation consists of $2^{n-1}$ positive rationals
$x_{n}^{(i)}$, $1\leq i\leq 2^{n-1}$. We denote this tree by $\mathcal{T}$, and
its $n$th generation by $\mathcal{T}^{(n)}$. The limitation of this tree to the
interval $[0,1]$ is the well-known Farey tree (albeit with different order).
Reading the tree line by line, this enumeration of $\Q_{+}$ starts with
\begin{eqnarray*}
{1\over 1}\,,\ {1\over 2}\,,\ {2\over 1}\,,\ {1\over 3}\,,\
{3\over 2}\,,\ {2\over 3}\,,\ {3\over 1}\,,\ {1\over 4}\,,\
{4\over 3}\,,\ {3\over 5}\,,\ {5\over 2}\,,\ {2\over 5}\,,\
{5\over 3}\,,\ {3\over 4}\,,\ {4\over 1}\,,\ \ldots
\end{eqnarray*}
This sequence was already investigated by Stern \cite{stern}, where one
encounters the definition of the Stern-Brocot tree. The sequence satisfies the
remarkable iteration discovered by M. Newman \cite{newman}:
$$
x_1=1,\quad x_{n+1}=1/(2[x_n]+1-x_n),
$$
thus giving an example of a simple recurrence which produces all positive
rationals, and answering affirmatively to a question by D.E. Knuth. The $n$th
generation of $\mathcal{T}$ consists of exactly those rational numbers whose
elements of the continued fraction sum up to $n$; this observation is due to
Stern \cite{stern}. Indeed, this can be easily inherited directly from the
definition. First, if rational number $\frac{a}{b}$ is represented as a
continued fraction $[a_{0},a_{1},...,a_{r}]$, then the map
$\frac{a}{b}\rightarrow\frac{a+b}{b}$ maps $\frac{a}{b}$ to
$[a_{0}+1,a_{1}...,a_{r}]$. Second, the map
$\frac{a}{b}\rightarrow\frac{a}{a+b}$ maps $\frac{a}{b}$ to
$[0,a_{1}+1,...,a_{r}]$ in case $\frac{a}{b}<1$, and to
$[1,a_{0},a_{1},...,a_{r}]$ in case $\frac{a}{b}>1$. This simple fact is of
utmost importance in our work: though it is not used in explicit form, this
highly motivates the investigations of moments $M_{L}$ and $m_{L}$, given by
(\ref{cont}). The sequence of numerators of the Calkin-Wilf tree
\begin{eqnarray*}
0,1,1,2,1,3,2,3,1,4,3,5,2,5,3,4,1,...
\end{eqnarray*}
is called the Stern diatomic sequence \cite{stern}, \cite{lehmer} and it
satisfies the recurrence relations
\begin{eqnarray}
s(0)=0,\quad s(1)=1,\quad s(2n)=s(n),\quad
s(2n+1)=s(n)+s(n+1).\label{stern}
\end{eqnarray}
\indent In the next section we will show that each generation of the
Calkin-Wilf tree possesses a distribution function $F_{n}(x)$, and $F_{n}(x)$
converges uniformly to $F(x)$. This is by far not a new fact. Nevertheless, we
include the short proof of it for the sake of self-containedness. The function
$F(x)$, as a distribution function, is uniquely determined by the functional
equation (\ref{distr}). This implies the explicit expression (\ref{min}) and the symmetry property $F(x)+F(1/x)=1$. The mean value of $F(x)$ was investigated by
several authors (\cite{steuding},\cite{wirsing2}) and was proved to be $3/2$.
We will obtain this result using quite a different method.
\subsection{Motivation and description of results} The aim of this paper is to give a different treatment of
Minkowski's $?(x)$. All papers so far were concerned with $F(x)$ as the function
{\it per se. }Nevertheless, it appears that there exist several natural
integral transforms of $F(x)$, which are analytic functions and which encode
certain substantial (in fact, all) information about the question mark function. Each of these
transforms is characterized by regularity properties and a functional equation.
Lastly and most importantly, let us point out that, surprisingly, there are
striking similarities and analogies between the results proved here as well as
in \cite{alkauskas1}, \cite{alkauskas3} with Lewis'-Zagier's \cite{lewis_zagier}
results on period functions for Maass wave forms. That work is an expanded and
clarified exposition of an earlier paper by Lewis \cite{lewis}. The concise
exposition of these objects, their properties and relations to Selberg zeta
function can be found in \cite{zagier}. The reader is strongly urged to compare
results in this work with those in \cite{lewis_zagier}. Let, for example,
$u(z)$ be a Maass wave form for ${\sf PSL}_{2}(\mathbb{Z})$ with spectral
parameter $s$. The similarity arises due to the fact that the limit value of
$u(z)$ on the real line, given by $u(x+iy)\sim y^{1-s}U(x)+y^{s}U(x)$ as
$y\rightarrow 0+$, satisfies (formal) functional equations $U(x+1)=U(x)$ and
$|x|^{2s-2}U\big{(}-\frac{1}{x}\big{)}=U(x)$. Thus, these are completely
analogous to the functional equations for $F(x)$, save the fact that $U(x)$ is
only a formal function - it is a distribution (e.g. a continuous functional in
properly defined space of functions). Thus, our objects $G(z)$,
$\mathfrak{m}(t)$ and $M(t)$ are analogues of objects $\psi(z)$, $g(w)$ and
$\phi(w)$ respectively (see Section 2 of this work and \cite{lewis_zagier}). In
\cite{alkauskas1} it is shown that in fact $L-$functions attached to Maass wave forms also do have an analogue in
context of the Minkowski question mark function. \\
\indent This work is organized as follows. In Section 2 we demonstrate some
elementary properties of the distribution function $F(x)$. Since the existence
of all moments is guaranteed by the exponential decay of the tail, our main
object is the generating function of moments, denoted by $G(z)$. In Section 3
we prove two functional equations for $G(z)$. In Section 4 we demonstrate the
uniqueness of solution of this functional equation (subject to regularity
conditions). Surprisingly, the Eisenstein series $G_{1}(z)$ appears on the
stage. In Section 5 we prove the integral equation for the exponential
generating function. In Section 6 a new class of functions emerging from
eigenfunctions of a Hilbert-Schmidt operator (we call them ``dyadic
eigenfunctions'') is introduced. These are dyadic analogues of functions
discovered by Wirsing \cite{wirsing} in connection with the
Gauss-Kuzmin-L\'{e}vy problem. In Section 7 we describe $p-$adic
distribution of rationals in the Calkin-Wilf tree. In the final section some
concluding remarks are presented.
\section{Some properties of the distribution}
The following Proposition was proved by many authors in various forms,
concerning (very related) Stern-Brocot, Farey or Calkin-Wilf trees, and this
seems to be a well-known fact about a distribution of rationals in these trees.
For the sake of completeness we present a short proof, since the functional
equations for $G(z)$ and $\mathfrak{m}(t)$ (see Sections 3 and 5) heavily
depend on the functional equation for $F(x)$ and are in fact reformulations of
these in different terms.
\begin{prop}\label{prop1}
Let $F_{n}(x)$ denote the distribution function of the $n$th generation, i.e.,
$$
F_{n}(x)=2^{1-n}\#\{j\,:\,x_{j}^{(n)}\leq x\}.
$$
Then uniformly $F_{n}(x)\rightarrow F(x)$.
Thus, $F(0)=0$, $F(\infty)=1$. Moreover, $F(x)$ is continuous, monotone and singular, i.e., $F'(x)=0$ almost everywhere.
\end{prop}

\noindent {\it Proof. }Let $x\geq1$. One half of the fractions in the $n+1-$st
generation do not exceed $1$, and hence also do not exceed $x$. Further,
$$
\frac{a+b}{b}\leq x\quad \iff\quad \frac{a}{b}\leq x-1.
$$
Hence,
$$
2F_{n+1}(x)=F_{n}(x-1)+1,\quad n\geq 1.
$$
Now assume $0<x<1$. Then
$$
\frac{a}{a+b}\leq x\quad\iff\quad \frac{a}{b}\leq\frac{x}{1-x}.
$$
Therefore,
$$
2F_{n+1}(x)=F_{n}\Big{(}\frac{x}{1-x}\Big{)}.
$$
The distribution function $F$, defined in the formulation of the Theorem,
satisfies the functional equation
\begin{eqnarray}
2F(x)=\left\{\begin{array}{c@{\qquad}l} F(x-1)+1 & \mbox{if}\quad x\geq 1,
\\ F({x\over 1-x}) & \mbox{if}\quad 0<x<1. \end{array}\right.\label{distr}
\end{eqnarray}
For instance, the second identity is equivalent to $2F(\frac{t}{t+1})=F(t)$ for all positive $t$.
If $t=[b_{0},b_{1},...]$, then $\frac{t}{t+1}=[0,1,b_{0},b_{1},..]$ for $t\geq 1$, and
 $\frac{t}{t+1}=[0,b_{1}+1,b_{2},...]$ for $t<1$, and the statement follows immediately.
\par

Now define $\delta_{n}(x)=F(x)-F_{n}(x)$. In order to prove the uniform convergence $F_n\to F$, it is sufficient to show that
\begin{equation}\label{vienas}
\sup_{x\geq 0}|\delta_{n}(x)|\leq 2^{-n}.
\end{equation}
It is easy to see that the assertion is true for $n=1$. Now suppose the estimate is true for $n$.
In view of the functional equation for both $F_n(x)$ and $F(x)$, we have
$$
2\delta_{n+1}(x)=\delta_{n}\Big{(}\frac{x}{1-x}\Big{)}
$$
for $0<x<1$, which gives $\sup_{0\leq x<1}|\delta_{n+1}(x)|\leq 2^{-n-1}$. Moreover, we have
$$
2\delta_{n+1}(x)=\delta_{n}(x-1)
$$
for $x\geq1$, which yields the same bound for $\delta_n(x)$ in the range $x\geq
1$. This proves (\ref{vienas}). As it was noted, the singularity of $F(x)$ was
proved in \cite{denjoy1} and it follows
from Khinchin results on metric properties of continued fractions. $\blacksquare$\\

Since $F(x)$ has a tail of exponential decay ($1-F(x)=O(2^{-x})$, as it is
clear from (\ref{min})), all moments do exist. Let
\[
M_{L}=\int\limits_{0}^{\infty}x^{L}\d F(x), \quad
m_{L}=\int\limits_{0}^{\infty}\Big{(}\frac{x}{x+1}\Big{)}^{L}\d F(x)=2\int\limits_{0}^{1}x^{L}\d F(x)=\int\limits_{0}^{1}x^{L}\d?(x).
\]
Therefore, $M_{L}$ and $m_{L}$ can also be defined as
\begin{eqnarray}
M_{L}=\lim_{n\rightarrow\infty}2^{1-n}\sum\limits_{a_{0}+a_{1}+...+a_{s}=n}[a_{0},a_{1},..,a_{s}]^{L},\quad
m_{L}=\lim_{n\rightarrow\infty}2^{2-n}\sum\limits_{a_{1}+...+a_{s}=n}[0,a_{1},..,a_{s}]^{L},\label{cont}
\end{eqnarray}
where the summation takes place over all rationals, whose elements of the
continued fraction sum up to $n$. These expressions highly motivate our
investigation of moments. Though the authors in \cite{grabner_tichy} considered
the moments of $F_{2}(x)$ (see the introduction), it is surprising that the
moments of Minkowski question mark function itself were never investigated.
Numerically, one has
\begin{align*}
M_{1}&=1.5,&\quad M_{2}&=4.290926,&\quad M_{3}&=18.556,&\quad
M_{4}&=107.03;\\
m_{1}&=0.5,&\quad m_{2}&=0.290926,&\quad m_{3}&=0.186389,&\quad
m_{4}&=0.126992.
\end{align*}
We will see that the generating function of $m_{L}$ possesses certain
fascinating properties. Let $\omega(x)$ be a continuous function of at most
polynomial growth: $\omega(x)=O(x^{T})$, $x\rightarrow\infty$. The functional
equation (\ref{distr}) gives $F(x+n)=1-2^{-n}+2^{-n}F(x)$, $x\geq0$. Hence
\[
\int\limits_{0}^{\infty}\omega(x)\d F(x)=
\sum\limits_{n=0}^{\infty}\int\limits_{0}^{1}\omega(x+n)\d F(x+n)=
\int\limits_{0}^{1}\sum\limits_{n=0}^{\infty}\frac{\omega(x+n)}{2^{n}}\d
F(x).
\]
Since, as noted above, $F(x)$ has a tail of exponential decay, this integral
does exist. Let $x=\frac{t}{t+1}$, $t\geq0$. Since
$F(\frac{t}{t+1})=\frac{1}{2}F(t)$, this change of variables gives
\[
\int\limits_{0}^{\infty}\omega(x)\d F(x)=
\sum\limits_{n=0}^{\infty}
\int\limits_{0}^{\infty}\frac{\omega(\frac{t}{t+1}+n)}{2^{n+1}}\d
F(t)
\]
(All changes of order of summation and integration are easily
justifiable minding the condition on $\omega(x)$). Let
$\omega(x)=x^{L}$, $L\in\mathbb{N}_{0}$. Then, if we denote
$B_{s}=\sum_{n=0}^{\infty}\frac{n^{s}}{2^{n+1}}$, we have
\[
\int\limits_{0}^{\infty}x^{L}\d F(x)=\int\limits_{0}^{\infty}
\sum\limits_{i=0}^{L}\Big{(}\frac{x}{x+1}\Big{)}^{i}
\binom{L}{i}B_{L-i}\d F(x).
\]
Whence the relation
\begin{eqnarray}
M_{L}=\sum\limits_{i=0}^{L}m_{i}\binom{L}{i}B_{L-i},\quad ,
\label{Mm} L\geq0.\label{Mm}
\end{eqnarray}
The exponential generating function of $B_{L}$ is
\[
B(t)=\sum\limits_{L=0}^{\infty}\frac{B_{L}}{L!}t^{L}=
\sum\limits_{L=0}^{\infty}\sum\limits_{n=0}^{\infty}
\frac{n^{L}t^{L}}{2^{N+1}L!}=
\sum\limits_{n=0}^{\infty}\frac{e^{nt}}{2^{n+1}}=\frac{1}{2-e^{t}}.
\]
Denote by $M(t)$ and $\mathfrak{m}(t)$ the corresponding exponential generating functions of the coefficients $M_{L}$ and $m_{L}$ respectively.
Accordingly,
\begin{eqnarray*}
M(t)=\int\limits_{0}^{\infty}e^{xt}\d F(x),\quad \mathfrak{m}(t)=
\int\limits_{0}^{\infty}\exp\Big{(}\frac{xt}{x+1}\Big{)}\d F(x)=2\int\limits_{0}^{1}e^{xt}\d F(x).
\end{eqnarray*}
The relation ($\ref{Mm}$) in terms of $M(t)$ and $\mathfrak{m}(t)$  reads as
\begin{eqnarray}
M(t)=\sum\limits_{L=0}^{\infty}\frac{M_{L}}{L!}t^{L}=\frac{1}{2-e^{t}}
\sum\limits_{L=0}^{\infty}\frac{m_{L}}{L!}t^{L}=\frac{1}{2-e^{t}}\mathfrak{m}(t).\label{exp}
\end{eqnarray}
We see that the function $\mathfrak{m}(t)$ is entire and $M(t)$ has a positive
radius of convergence. The last identity implies the asymptotic formula for
$M_{L}$.

\begin{prop} For $L\in\N_0$,
\begin{eqnarray*}
M_{L}&=&\frac{\mathfrak{m}(\log 2)}{2\log 2}\Big{(}\frac{1}
{\log2}\Big{)}^{L}L!+O_{\varepsilon}\Big{(}((4\pi^{2}+(\log 2)^{1/2}-\varepsilon)^{-L}\Big{)}L!\\
&=& \Big{(}\frac{\mathfrak{m}(\log 2)}{2\log 2}\Big{(}\frac{1}
{\log2}\Big{)}^{L}+O(6.3^{-L})\Big{)}L!
\end{eqnarray*}
\label{prop2}
\end{prop}

\noindent {\it Proof. }By Cauchy's formula, for any sufficiently small $r$,
$$
M_L={L!\over 2\pi i}\int_{\vert z\vert=r}{M(z)\over z^{L+1}}\d z.
$$
Changing the path of integration, we get by the calculus of residues
$$
M_L=-\mbox{Res}_{z=\log 2}\left({m(z)\over (2-e^z)z^{L+1}}\right)-{L!\over 2\pi i}\int_{\vert z\vert=R}{m(z)\over 2-e^z}{\d z\over z^{L+1}},
$$
where $R$ satisfies $\log 2<R<\vert \log 2+2\pi i\vert$ (which means that there
is exactly one simple pole of the integrand located in the interior of the
circle $\vert z\vert=R$). It is easily seen that the residue coincides with the
main term in the formula of the Lemma; the error term follows from estimating
the integral. $\blacksquare$
\par\bigskip
Also, (\ref{exp}) gives the inverse to linear equations ($\ref{Mm})$:
\begin{eqnarray}
m_{L}=M_{L}-\sum\limits_{s=0}^{L-1}M_{s}\binom{L}{s}, \quad
L\geq0.\label{inv}
\end{eqnarray}
Since $B(t)(2-e^{t})=1$, the coefficients $B_{L}$ can be calculated
recursively: $B_{L}=\sum_{s=0}^{L-1}\binom{L}{s}B_{s}$. Thus, $B_{0}=1$,
$B_{1}=1$, $B_{2}=3$, $B_{3}=13$, $B_{4}=75$, $B_{5}=541$. This sequence has
number $A000670$ in \cite{sloane}, and traces
its history back from Cayley.\\

In the future we will consider integrals which involve $\mathfrak{m}(t)$, and
hence we need the evaluation of this function for negative $t$.
\begin{lem} Let $C=e^{-\sqrt{\log 2}}=0.4349_{+}$. Then
$C^{2\sqrt{t}}\ll \mathfrak{m}(-t)\ll C^{\sqrt{t}}$ as $t\rightarrow\infty$.
\label{lem1}
\end{lem}
{\it Proof. }In fact,
$\mathfrak{m}(-t)=\int_{0}^{\infty}\exp(-\frac{xt}{x+1})\d F(x)$. Hence,
$\mathfrak{m}(t)$ is positive for $t\in\mathbb{R}$. Let $0<M<1$. Since
$1-F(x)\asymp2^{-x}$ as $x\rightarrow\infty$, and $F(x)+F(1/x)=1$,
\[
\mathfrak{m}(-t)=\Big{(}\int\limits_{0}^{M}+\int\limits_{M}^{\infty}\Big{)}\exp(-\frac{xt}{x+1})\d
F(x)\ll 2^{-1/M}+\exp\Big{(}-\frac{Mt}{M+1}\Big{)}.
\]
This is valid for every $M<1$ and a universal constant. A choice
$M=\frac{\sqrt{\log 2}}{\sqrt{t}}$ gives the desired upper bound. To obtain
the lower bound, note that
\begin{eqnarray*}
\mathfrak{m}(-t)>\int\limits_{0}^{M}\exp(-\frac{xt}{x+1})\d
F(x)\gg 2^{-1/M}\cdot\exp\Big{(}-\frac{Mt}{M+1}\Big{)}.
\end{eqnarray*}
The same choice for $M$ establishes the lower bound. Naturally, similar
evaluation holds for the derivative, since
$\mathfrak{m}'(-t)=\int_{0}^{\infty}\frac{x}{x+1}\exp(-\frac{xt}{x+1})\d F(x)$. $\blacksquare$\\

We will prove one property of the function $\mathfrak{m}(t)$ which represents
the symmetry of $F$ given by $F(x)+F(1/x)=1$.
\begin{prop}
We have: $\mathfrak{m}(t)=e^{t}\mathfrak{m}(-t)$. \label{prop3}
\end{prop}
{\it Proof. }In fact,
\[
\mathfrak{m}(t)=\int\limits_{0}^{\infty}\exp\Big{(}\frac{xt}{x+1}\Big{)}\d
F(x)=
-\int\limits_{0}^{\infty}\exp\Big{(}\frac{t/x}{1/x+1}\Big{)}\d
F(1/x)\]
\[=\int\limits_{0}^{\infty}\exp\Big{(}\frac{t}{x+1}\Big{)}\d F(x)
=e^{t}\int\limits_{0}^{\infty}\exp\Big{(}-\frac{xt}{x+1}\Big{)}\d
F(x)=\mathfrak{m}(-t)e^{t}.\quad\blacksquare
\]
Whence the relations
\begin{eqnarray*}
m_{L}=\sum\limits_{s=0}^{L}\binom{L}{s}(-1)^{s}m_{s},
\quad, L\geq0.
\end{eqnarray*}
Thus, $m_{1}=m_{0}-m_{1}$, which gives $m_{1}=1/2$, and this implies
$M_{1}=3/2$. Also, $2m_{3}=-1/2+3m_{2}$. These linear relations are further
investigated in \cite{alkauskas1}.
\section{The dyadic period function $G(z)$}

We introduce the generating power function of moments
\begin{eqnarray*}
\mathcal{M}(z)=\sum\limits_{L=0}^{\infty}m_{L}z^{L}.
\end{eqnarray*}
{\it A priori, }this series converges in the unit circle. Recall that
$\int_{0}^{\infty}x^{n}e^{-x}\d x=\Gamma(n+1)=n!$. Thus, for real $z<1$, the
symmetry relation for $\mathfrak{m}(t)$ gives:
\begin{eqnarray}
\mathcal{M}(z)=\int\limits_{0}^{\infty}\mathfrak{m}(zt)e^{-t}\d t=
\int\limits_{0}^{\infty}m(-zt)e^{-t(1-z)}\d t=\nonumber\\
\int\limits_{0}^{\infty}m\Big{(}t\frac{z}{z-1}\Big{)}
\frac{1}{(1-z)}e^{-t}\d t=
\mathcal{M}\Big{(}\frac{z}{z-1}\Big{)}\frac{1}{1-z}.\label{symm}
\end{eqnarray}
Both integrals converge for $z<1$ (since $m_{L}\leq 1$, $|m(z)|\leq e^{z}$),
hence for these values of $z$ we have the above identity. The function
$\mathcal{M}(z)$ was initially defined for $|z|<1$; nevertheless, the above
identity gives us holomorphic continuation of $\mathcal{M}(z)$ to the half
plane $\Re{z}<1/2$.
\begin{lem} The function $\mathcal{M}(z)$ can be analytically
continued to the domain $\mathbb{C}\backslash\mathbb{R}_{x> 1}$. \label{lem2}
\end{lem}
{\it Proof. }In fact, $\mathfrak{m}(t)=\int_{0}^{\infty}\exp(\frac{x}{x+1}t)\d
F(x)$. As noted above, $|\mathfrak{m}(t)|\leq e^{t}$ for positive $t$
(actually, Lemma \ref{lem1} combined with Proposition \ref{prop3} gives a slightly better estimate). Therefore, for
real $z<1$ we have:
\[
\mathcal{M}(z)=\int\limits_{0}^{\infty}\int\limits_{0}^{\infty}
\exp\Big{(}\frac{x}{x+1}zt\Big{)}e^{-t}\d F(x)\d t=
\int\limits_{0}^{\infty}\frac{1}{1-\frac{x}{x+1}z}\d F(x).
\]
We already obtained the analytic continuation of $\mathcal{M}$ to the region
$\{|z|<1\}\bigcup\{\Re z<1/2\}$. Let $z=\sigma+iy$ with $y>0$ and
$\sigma\geq1/2$. In the small neighborhood of $z$ the imaginary part is
bounded: $y\geq y_{0}>0$, and also the real part is bounded:
$\sigma\leq\sigma_{0}$. In this neighborhood the integral converges uniformly,
because we have the estimate
\begin{eqnarray*}
\Big{|}(1-\frac{x}{x+1}z)\Big{|}\geq
\max\{\Big{|}1-\frac{x}{x+1}\sigma\Big{|},
\Big{|}\frac{x}{x+1}t\Big{|}\}.
\end{eqnarray*}
For $0\leq x\leq 1/\sigma$ this gives the bound $\frac{1}{\sigma_{0}+1}$, and
for $x>1/\sigma$ this gives the bound $\frac{t_{0}}{\sigma_{0}+1}$. Hence, in
this neighborhood the function under the integral is uniformly bounded, which
proves the uniform convergence of the integral and the statement of Lemma. $\blacksquare$\\

The system ($\ref{Mm}$) gives us the expression of $M_{L}$ in terms of $m_{s}$.
In fact, there exists one more system which is independent of the distribution
$F(x)$; it simply encodes the relation among functions
$\big{(}\frac{x}{1-x}\big{)}^{L}$ and $x^{s}$, given by
\[
\Big{(}\frac{x}{1-x}\Big{)}^{L}=\sum\limits_{s\geq
L}\binom{s-1}{L-1}x^{s},\quad L\geq 1,\quad 0\leq x<1.
\]
A change $x=\frac{t}{t+1}$ gives
\[
t^{L}=\sum\limits_{s\geq
L}\binom{s-1}{L-1}\Big{(}\frac{t}{t+1}\Big{)}^{s}
\quad L\geq 1, \quad t\geq 0.
\]
And ultimately,
\begin{eqnarray}
M_{L}=\sum\limits_{s\geq L}\binom{s-1}{L-1}m_{s}.\label{rys}
\end{eqnarray}
For the convenience, we introduce a function
\begin{eqnarray}
G(z)=\frac{\mathcal{M}(z)-1}{z}=\sum_{L=1}^{\infty}m_{L}z^{L-1}=
\int\limits_{0}^{\infty}\frac{\frac{x}{x+1}}{1-\frac{x}{x+1}z}\d F(x)=
2\int\limits_{0}^{1}\frac{x}{1-xz}\d F(x).\label{rep}
\end{eqnarray}
\indent Next theorem is our main result about $G(z)$. The power series
converges in the disc $|z|\leq 1$ (including the boundary, as can be inherited
from (\ref{rys}); moreover, this implies that there exist all left derivatives
of $G(z)$ at $z=1$). The integral converges in the cut plane
$\mathbb{C}\setminus(1,\infty)$.
\begin{thm}\label{thm1}
Let $m_{L}=\int_{0}^{\infty}(\frac{x}{x+1})^{L}\d F(x)$. Then the
generating power function, defined as
$G(z)=\sum_{L=1}^{\infty}m_{L}z^{L-1}$, has an analytic
continuation to the domain $\mathbb{C}\backslash\mathbb{R}_{x>
1}$. It satisfies the functional equation
\begin{eqnarray}
-\frac{1}{1-z}-\frac{1}{(1-z)^{2}}G\Big{(}\frac{1}{1-z}\Big{)}+2G(z+1)=G(z),
\label{ggg}
\end{eqnarray}
and also
the symmetry property
$$
G(z+1)=-\frac{1}{z^{2}}G\Big{(}\frac{1}{z}+1\Big{)}-\frac{1}{z}.
$$
Moreover, $G(z)=o(1)$ as $z\rightarrow\infty$
and the distance from $z$ to $\mathbb{R}_{+}$ tends to infinity.
\end{thm}
{\it Proof. }In analogy to $\mathcal{M}(z)$, for real $z<0$ define the
following function: $\mathcal{M}_{0}(z)=\int_{0}^{\infty}M(zt)e^{-t}\d t$. In
view of $(\ref{exp})$, this integral converges for real $z<0$. Thus,
\[
\mathcal{M}_{0}(z)=\int\limits_{0}^{\infty}\int\limits_{0}^{\infty}
\exp(xzt)e^{-t}\d F(x)\d t=
\int\limits_{0}^{\infty}\frac{1}{1-xz}\d F(x).
\]
From argument akin to the one used in proving Lemma \ref{lem2} we deduce that $\mathcal{M}_{0}(z)$
extends as an analytic function to the region
$\mathbb{C}\backslash\mathbb{R}_{> 0}$. In this domain we see that
\begin{eqnarray}
\frac{\mathcal{M}_{0}(z)-1}{z}=\frac{\mathcal{M}(z+1)-1}{z+1},\label{om}
\end{eqnarray}
which is the consequence of an algebraic identity
\[
\Big{(}\frac{1}{1-xz}-1\Big{)}\cdot\frac{1}{z}=
\Big{(}\frac{1}{1-\frac{x}{x+1}(z+1)}-1\Big{)}\cdot\frac{1}{z+1}.
\]
The relation $(\ref{om})$ is independent of the specific distribution function,
it simply encodes the information contained in $(\ref{rys})$ about the relation
of $x^{L}$ to $(x/(x+1))^{l}$. On the other hand, the specific information
about $F(x)$ is encoded in $(\ref{Mm})$ or $(\ref{exp})$. The comparison of
these two relations gives the desired functional equation for $G(z)$. In fact,
for real $t<0$ the following estimate follows from ($\ref{exp}$) and Lemma
\ref{lem1}: $|M(t)|=|\mathfrak{m}(t)(2-e^{t})^{-1}|\leq|\mathfrak{m}(t)|\ll 1$;
and thus for real $z<0$ we have:
\[
\mathcal{M}(z)=\int\limits_{0}^{\infty}\mathfrak{m}(zt)e^{-t}\d t=
\int\limits_{0}^{\infty}(2-e^{zt})M(zt)e^{-t}\d t=\]
\[
2\mathcal{M}_{0}(z)-\int\limits_{0}^{\infty}M(zt)e^{-t(1-z)}\d t=
2\mathcal{M}_{0}(z)-\mathcal{M}_{0}\Big{(}\frac{z}{1-z}\Big{)}\frac{1}{1-z}.
\]
Finally,  the substitution $(\ref{om})$ gives us the functional equation
\begin{eqnarray*}
\frac{1-z}{1+z}-\frac{z}{1-z}\mathcal{M}\Big{(}\frac{1}{1-z}\Big{)}+
2\frac{z}{z+1}\mathcal{M}(z+1)=\mathcal{M}(z).
\end{eqnarray*}

The principle of analytic continuation implies that this equation  should be
satisfied for all values of arguments in the region of holomorphicity of
$\mathcal{M}(z)$. Direct inspection shows that for
$G(z)=\frac{\mathcal{M}(z)-1}{z}$ this equation reads as (\ref{ggg}). Also, the
symmetry property is a reformulation or (\ref{symm}). This proves the first
part of
the Theorem.\\
\indent Obviously, the last assertion follows from the
integral representation of $G(z)$ given by (\ref{rep}). $\blacksquare$\\

We call $G(z)$ {\it the dyadic period function}, since its
functional equation is completely parallel to a three term functional equations
which are satisfied by rational period functions and period functions
associated with Maass wave forms \cite{lewis_zagier}. The word ``dyadic''
refers to the binary origin of the distribution function $F(x)$. Indeed,
thorough inspection shows that the multiplier $2$ in equations (\ref{eigen})
and (\ref{ggg}) emerges exactly from the fact that every generation of
$\mathcal{T}$ has twice as many members as a previous generation.\\
\section{Uniqueness of $G(z)$}
In this section we prove the uniqueness of a function which have properties
described in Theorem \ref{thm1}. Note that two functional equations for $G(z)$
can be merged into a single one. It is easy to check that
\begin{eqnarray}
\frac{1}{z}+\frac{1}{z^{2}}G\Big{(}\frac{1}{z}\Big{)}+2G(z+1)=G(z)\label{sim}
\end{eqnarray}
is equivalent to both together. In fact, the change $z\mapsto 1/z$ in the
equation (\ref{sim}) gives the symmetry property, and application of it to the term $G(1/z)$ of the above recovers the functional equation (\ref{ggg}).
\begin{prop}
The function which satisfies the conditions of Theorem \ref{thm1} is unique.
\label{prop4}
\end{prop}
{\it Proof. }Suppose there exist two such functions. Then their difference
$G_{0}(z)$ has the same behavior at infinity, and it satisfies the homogenic
form of the equation ($\ref{sim}$). Let $M=\sup_{-1\leq x\leq
0}|G_{0}(x)|=|G_{0}(x_{0})|$, $x_{0}\in[-1,0]$. We will show that $M=0$; by the
principle of analytic continuation this will imply that $G_{0}(z)\equiv 0$. Let
$z$ be real, $-1\leq z\leq 0$. Let us substitute $z\mapsto z-n$ in the equation
($\ref{sim}$), $n\in\mathbb{N}$, $n\geq 1$, and divide it by $2^{n}$. Thus, we
obtain:
\begin{eqnarray}
\frac{G_{0}(z-n)}{2^{n}}-\frac{G_{0}(z-n+1)}{2^{n-1}}=
\frac{1}{2^{n}(z-n)^{2}}G_{0}\Big{(}\frac{1}{z-n}\Big{)}. \label{akagkw}
\end{eqnarray}
Note that for $z$ in the interval $[-1,0]$, $\frac{1}{z-n}$ belongs to the same
interval as well. Now sum this over $n\geq 1$. The series on both sides are
absolutely convergent, minding the behavior of $G_{0}(z)$ at infinity.
Therefore,
\begin{eqnarray*}
-G_{0}(z)=\sum_{n=1}^{\infty}\frac{1}{2^{n}(z-n)^{2}}
G_{0}\Big{(}\frac{1}{z-n}\Big{)}.
\end{eqnarray*}
The evaluation of the right hand side gives:
\begin{eqnarray*}
|G_{0}(z)|\leq\sum_{n=1}^{\infty}\frac{1}{2^{n}n^{2}}M=
\Big{(}\frac{\pi^{2}}{12}-\frac{1}{2}\log^{2}2\Big{)}M \quad
\text{for }-1\leq z\leq 0.
\end{eqnarray*}
The constant is $<1$. Thus, unless $M=0$, this is contradictory for $z=x_{0}$.
This proves the
Proposition. $\blacksquare$\\

Note the similarity between (\ref{akagkw}) and the expression for the
Gauss-Kuzmin-Wirsing operator $\mathbf{W}$. The latter is defined for bounded
smooth functions $f:[0,1]\rightarrow\mathbb{R}$ by the formula
\[
[\mathbf{W}f](x)=\sum\limits_{k=1}^{\infty}\frac{1}{(k+x)^{2}}
f\Big{(}\frac{1}{k+x}\Big{)}.
\]
The eigenvalue $1$ corresponds to the function $\frac{1}{1+x}$ (see
\cite{khinchin} chapter III, for Kuzmin's treatment). The second largest
eigenvalue $-0.303663...$ (the Wirsing constant) leads to a function with no
analytic expression known \cite{wirsing}; this eigenvalue determines the speed
of convergence of iterates $[\mathbf{W}^{(n)}f](x)$ to $\frac{c}{1+x}$ (for
certain $c\in\mathbb{R}$). The spectral analysis of our operator is presented
in Section 6. Paper \cite{alkauskas1}
contains much more details and results in this direction.\\

Let $\Im z>0$. We remind that the Eisenstein series of weight $2$ for ${\sf PSL}_{2}(\mathbb{Z})$ is defined as \cite{serre}
\[
G_{1}(z)=\sum\limits_{n\in\mathbb{Z}}{\sum\limits_{m\in\mathbb{Z}}}^{'}\frac{1}{(m+nz)^{2}};
\]
(mind the order of summation, since the series is not absolutely
convergent). This series has the following Fourier expansion:
if $q=e^{2\pi iz}$, then
\[
G_{1}(z)=\frac{\pi^{2}}{3}-8\pi^{2}\sum\limits_{n=1}^{\infty}\sigma_{1}(n)q^{n}.
\]Then this function is not completely modular, but we
have the following identities (\cite{serre}, chapter VII):
\[
G_{1}(z+1)=G_{1}(z),\quad G_{1}(-1/z)=z^{2}G_{1}(z)-2\pi i z.
\]
Note that for $\Im z>0$, all arguments in (\ref{ggg})
simultaneously belong to the upper half plane. It is surprising
(but not coincidental) that the function $\frac{i}{2\pi}G_{1}(z)$
satisfies the functional equation (\ref{ggg}) for $\Im z>0$ (see
the remarks in Section 8 about possible connections in idelic
setting). To check this statement, note that
\[
\frac{i}{2\pi} G_{1}\Big{(}-\frac{1}{z-1}\Big{)}=
\frac{i}{2\pi}\Big{(}(z-1)^{2}G_{1}(z-1)- 2\pi
i(z-1)\Big{)}=\frac{i}{2\pi}(1-z)^{2}G_{1}(z)-(1-z).
\]
Thus, plugging this into ($\ref{ggg}$), we obtain an identity. If
we define $G_{1}(z)=\overline{G_{1}(\overline{z})}$ for $\Im z<0$,
one checks directly that the symmetry property is also satisfied.
This is a surprising phenomena. See the last section of \cite{alkauskas1} for more speculations on this topic,
where the space of dyadic period functions in the upper half plane (denoted by ${\tt DPF}^{0}$) is introduced.\\

We end this section with presenting a system of linear equations which the moments $m_{L}$ do satisfy. This system is derived from the three therm functional equation
(\ref{sim}) and is a superior result in numerical calculations: whereas directly
from the definition we can recover only a few digits of the moments, this method
allows to calculate up to $60$ digits and more.
\begin{prop} Denote
$c_{L}=\sum_{n=1}^{\infty}\frac{1}{2^{n}n^{L}}=\text{Li}_{L}(\frac{1}{2})$. The
moments $m_{s}$ satisfy the infinite system of linear equations
\begin{eqnarray*}
m_{s}=\sum\limits_{L=0}^{\infty}(-1)^{L}c_{L+s}\binom{L+s-1}{s-1}m_{L},
\quad s\geq 1.
\end{eqnarray*}
\label{prop5}
\end{prop}
{\it Proof. }Indeed, for $\Re z\leq 0$ we have (recall that $m_{0}=1$):
\begin{eqnarray*}
-G(z)=\sum\limits_{n=1}^{\infty}\frac{1}{2^{n}(z-n)}+\sum_{n=1}^{\infty}\frac{1}{2^{n}(z-n)^{2}}
G\Big{(}\frac{1}{z-n}\Big{)}=\sum_{n=1}^{\infty}\frac{1}{2^{n}}
\sum\limits_{L=0}^{\infty}m_{L}\Big{(}\frac{1}{z-n}\Big{)}^{L+1}.
\end{eqnarray*}
This series is absolutely and uniformly convergent for $\Re z\leq 0$, as is
implied by (\ref{rys}). We obtain the needed result after
taking the $s$th left derivative at $z=0$. $\blacksquare$\\

Numerical calculations are presented in \cite{alkauskas3}. This method gives
high precision values for other constants, including the Kinney's constant.
\section{Exponential generating function $\mathfrak{m}(t)$}

The aim of this section is to interpret (\ref{sim}) in terms of
$\mathfrak{m}(t)$. The following Theorem, along with the boundary condition
$\mathfrak{m}(0)=1$ and regularity property as in Lemma \ref{lem1}, uniquely determines
the function $\mathfrak{m}(t)$ .
\begin{thm}
The function $\mathfrak{m}(s)$ satisfies the integral equation
\begin{eqnarray}
\mathfrak{m}(-s)=(2e^{s}-1)\int\limits_{0}^{\infty}\mathfrak{m}'(-t)J_{0}(2\sqrt{st})\d
t, \quad s\in\mathbb{R}_{+},\label{bes}
\end{eqnarray}
where $J_{0}(*)$ stands for the Bessel function:
$J_{0}(z)=\frac{1}{\pi}\int_{0}^{\pi}\cos(z\sin x)\d x$.

\label{thm2}
\end{thm}
{\it Proof. }For $\Re z<1$, we have that $G(z)=\int_{0}^{\infty}m'(zt)e^{-t}\d
t$. Thus,
\begin{eqnarray*}
G(z)=-\frac{1}{z}\int\limits_{0}^{\infty}\mathfrak{m}'(-t)e^{t/z}\d t \text{
for } \Re z<0,\quad G(z)=\frac{1}{z}\int\limits_{0}^{\infty}m'(t)e^{-t/z}\d t
\text{ for
 } 0<\Re z<1.
\end{eqnarray*}
Thus, the functional equation for $G(z)$ in the region $\Re z<-1$
in terms of $\mathfrak{m}'(t)$ reads as
\begin{eqnarray}
\frac{1}{z}=\int\limits_{0}^{\infty}\mathfrak{m}'(-t)
\Big{(}\frac{2}{z+1}e^{\frac{t}{z+1}}+\frac{1}{z}e^{tz}-
\frac{1}{z}e^{\frac{t}{z}}\Big{)}\d t.
\end{eqnarray}\label{inteq}
Now, multiply this by $e^{-sz}$ and integrate over $\Re z=-\sigma<-1$, where
$s>0$ is real. We have (\cite{lavr}, p. 465)
\[
\int\limits_{-\sigma-i\infty}^{-\sigma+i\infty}\frac{e^{-sz}}{z}\d
z =-2\pi i;
\]
\begin{eqnarray*}
2\int\limits_{-\sigma-i\infty}^{-\sigma+i\infty}\frac{e^{\frac{t}{z+1}-sz}}{z+1}\d
z=
-2e^{s}\int\limits_{\sigma-1-i\infty}^{\sigma-1+i\infty}\frac{e^{sz-\frac{t}{z}}}{z}\d
z= -2e^{s}\int\limits_{\sigma_{0}-i\infty}^{\sigma_{0}+i\infty}
\frac{e^{\sqrt{st}z-\frac{\sqrt{st}}{z}}}{z}\d z=-4\pi i
e^{s}J_{0}(2\sqrt{st}),
\end{eqnarray*}
where $\sigma_{0}=(\sigma-1)\sqrt{\frac{t}{s}}>0$, and
$J_{\lambda}(*)$ stands for the Bessel function (see \cite{lavr},
p. 597 for the representation of the Bessel function by this
integral). Further,
\begin{eqnarray*}
\int\limits_{-\sigma-i\infty}^{-\sigma+i\infty}\frac{e^{(t-s)z}}{z}\d
z=
\begin{cases}
-2\pi i & \text{ if }s>t,\\
0 & \text { if }s<t,\quad
\end{cases}
\int\limits_{-\sigma-i\infty}^{-\sigma+i\infty}\frac{e^{\frac{t}{z}-sz}}{z}\d
z= -2\pi iJ_{0}(2\sqrt{st}).
\end{eqnarray*}
Thus, eventually
\begin{eqnarray*}
-2\pi i=-2\pi i\int\limits_{0}^{\infty}\mathfrak{m}'(-t)
(2e^{s}-1)J_{0}(2\sqrt{st})\d t -2\pi i\int\limits_{0}^{s}\mathfrak{m}'(-t)\d
t;
\end{eqnarray*}
since $\mathfrak{m}(0)=1$, this proves the Proposition. $\blacksquare$\\

Thus, we have obtained an integral equation for $\mathfrak{m}(s)$, which
corresponds to the functional equation (\ref{sim}) for $G(z)$. They are both,
in fact, mere reformulations of (\ref{min}) in different terms.

\section{Dyadic eigenfunctions}
In this Section we introduce the sequence of functions $G_{\lambda}(z)$, which
satisfy the functional equation analogous to ($\ref{sim}$).\\

Since $J'_{0}(*)=-J_{1}(*)$, $J_{1}(0)=0$, integration by parts in (\ref{bes})
leads to
\begin{eqnarray}
\int\limits_{0}^{\infty}\frac{\mathfrak{m}(-t)}{\sqrt{t}} J_{1}(2\sqrt{st})\d
t=\frac{1}{\sqrt{s}}-\frac{\mathfrak{m}(-s)}{\sqrt{s}(2e^{s}-1)}.\label{hank}
\end{eqnarray}
Recall that the Hankel transform of degree $\nu>-1/2$ of the function $f(r)$
(provided that $\int_{0}^{\infty}f(r)\sqrt{r}\d r$ converges absolutely) is
defined as
\[
g(\rho)=\int\limits_{0}^{\infty}f(r)J_{\nu}(r\rho)r\d r,
\]
where $J_{\nu}(*)$ stands for the $\nu$th Bessel function. The inverse is given by the Hankel inversion formula with exactly the same kernel (\cite{watson}, chapter
XIV, section 14.4.). Thus, after a proper change of variables, Hankel transform
reads as
\begin{eqnarray*}
g(\rho)=\int\limits_{0}^{\infty}f(r)J_{\nu}(2\sqrt{r\rho})\d
r\Leftrightarrow
f(r)=\int\limits_{0}^{\infty}g(\rho)J_{\nu}(2\sqrt{r\rho})\d\rho.
\end{eqnarray*}
Thus, application of this inversion to the identity (\ref{hank}) yields
\[
\frac{\mathfrak{m}(-s)}{\sqrt{s}}=\int\limits_{0}^{\infty}\frac{J_{1}(2\sqrt{st})}{\sqrt{t}}\d
t-\int\limits_{0}^{\infty}\frac{\mathfrak{m}(-t)}{\sqrt{t}(2e^{t}-1)}
J_{1}(2\sqrt{st})\d t.
\]
The first integral on the r.h.s. is equal t
$-\frac{1}{\sqrt{s}}J_{0}(2\sqrt{st})\mid_{t=0}^{\infty}=\frac{1}{\sqrt{s}}$.
Let $\psi(s)=(2e^{s}-1)^{1/2}$. Then this equation can be rewritten as
\[
\frac{\mathfrak{m}(-s)}{\sqrt{s}\psi(s)}=\frac{1}{\sqrt{s}\psi(s)}-
\int\limits_{0}^{\infty}\frac{\mathfrak{m}(-t)}{\sqrt{t}\psi(t)}\cdot
\frac{J_{1}(2\sqrt{st})}{\psi(s)\psi(t)}\d t.
\]
Hence, if we denote
\[
\frac{J_{1}(2\sqrt{st})}{\psi(s)\psi(t)}=K(s,t), \quad
\frac{\mathfrak{m}(-s)-1}{\sqrt{s}\psi(s)}=\mathbf{m}(s),
\]
we obtain a second type Fredholm integral equation with symmetric kernel
(\cite{kolm}, chapter 9):
\begin{eqnarray*}
\mathbf{m}(s)=\ell(s)-\int\limits_{0}^{\infty}\mathbf{m}(t)K(s,t)\d
t,\\
\text{ where }\ell(s)=-\frac{1}{\psi(s)}\int\limits_{0}^{\infty}
\frac{J_{1}(2\sqrt{st})}{\sqrt{t}(2e^{t}-1)}\d t=\frac{1}{\sqrt{s}\psi(s)}\Big{(}\sum\limits_{n=1}^{\infty}e^{-s/n}2^{-n}-1\Big{)}.
\end{eqnarray*}
The behavior at infinity of the Bessel function is given by an
asymptotic formula
\[
J_{1}(x)\sim\Big{(}\frac{2}{\pi
x}\Big{)}^{1/2}\cos\Big{(}x-\frac{3}{4}\pi\Big{)}
\]
(\cite{watson}, chapter VII, section 7.1). Therefore, obviously,
\[
\int\limits_{0}^{\infty}\int\limits_{0}^{\infty}|K(s,t)|^{2}\d s\d
t<\infty,\quad \int\limits_{0}^{\infty}|\ell(s)|^{2}\d s<\infty.
\]
Thus, the operator associated with the kernel $K(s,t)$ is the Hilbert-Schmidt
operator (\cite{kolm}, p. 532). The theorem of Hilbert-Schmidt (\cite{kolm}, p.
283) states that the solution of this type of integral equations reduces to
finding the eigenvalues $\lambda$ and the eigenfunctions $A_{\lambda}(s)$. We
postpone the solution of this integral equation for the future. Till the end of
this section we deal only with eigenfunctions. The integral operator,
consequently, is a compact self-conjugate operator in the Hilbert space, it
possesses a complete orthogonal system of eigenfunctions $\lambda$, all $\lambda$ are
real and $\lambda_{n}\rightarrow0 $, as $n\rightarrow\infty$. If we denote
$A_{\lambda}(s)\psi(s)=B_{\lambda}(s)$, then the equation for an eigenfunction
reads as
\[
\int\limits_{0}^{\infty}B_{\lambda}(t)\frac{J_{1}(2\sqrt{st})}{2e^{t}-1}\d
t= \lambda B_{\lambda}(s).
\]
This gives $B_{\lambda}(0)=0$. Since $A_{\lambda}(s)\in
L_{2}(0,\infty)$, and $J_{1}(*)$ is bounded, this implies that
$B_{\lambda}(s)$ is uniformly bounded for $s\geq0$ as well.
Moreover, since the Taylor expansion of $J_{1}(*)$ contains only
odd powers of the variable, $B_{\lambda}(s)\sqrt{s}$ has a Taylor
expansion with center $0$ and is an entire function+. Now, multiply this by
$\sqrt{s}e^{-s/z}$, $z>0$, and integrate over
$s\in\mathbb{R}_{+}$. The Laplace transform of
$\sqrt{s}J_{1}(2\sqrt{s})$ is $\frac{1}{z^{2}}e^{-1/z}$
(\cite{lavr}, p. 503). Thus, we obtain
\begin{eqnarray}
\frac{1}{\lambda}\int\limits_{0}^{\infty}\frac{B_{\lambda}(t)\sqrt{t}}{2e^{t}-1}e^{-tz}\d
t=\frac{1}{z^{2}}\int\limits_{0}^{\infty}B_{\lambda}(s)\sqrt{s}e^{-\frac{s}{z}}\d
s.\label{hil}
\end{eqnarray}
Denote by $G_{\lambda}(-z)$ the function on both sides of this equality. Thus,
$G_{\lambda}(z)$ is defined at least for $\Re z\leq0$. Since
$2e^{t(z+1)}-e^{tz}=(2e^{t}-1)e^{tz}$, we have
\begin{eqnarray*}
\lambda\Big{(}2G_{\lambda}(z+1)-G_{\lambda}(z)\Big{)}=
\int\limits_{0}^{\infty}B_{\lambda}(t)\sqrt{t}e^{tz}\d
t=\frac{1}{z^{2}}G_{\lambda}(1/z).
\end{eqnarray*}
Therefore, we have proved the first part of the following Theorem.
\begin{thm}
For every eigenvalue $\lambda$ of the integral operator, associated with the
kernel $K(s,t)$, there exists at least one holomorphic function $G_{\lambda}(z)$
(defined for $z\in\mathbb{C}\setminus\mathbb{R}_{>1}$), such that the following
holds:
\begin{eqnarray}
2G_{\lambda}(z+1)=G_{\lambda}(z)+\frac{1}{\lambda
z^{2}}G_{\lambda}\Big{(}\frac{1}{z}\Big{)}.\label{eigen}
\end{eqnarray}
Moreover, $G_{\lambda}(z)$ for $\Re z<0$ satisfies all regularity
conditions imposed by it being an image under the Laplace transform
(\cite{lavr}, p. 469).\\
Conversely: for every $\lambda$, such that there exists a non-zero function, which
satisfies $(\ref{eigen})$ and these conditions, $\lambda$ is the eigenvalue of
this operator. The set of all possible $\lambda$'s is countable, and
$\lambda_{n}\rightarrow 0$, as $n\rightarrow\infty$. \label{thm3}
\end{thm}
{\it Proof. }The converse is straightforward, since, by the requirement,
$G_{\lambda}(z)$ for $\Re z\leq 0$ is a Laplace image of a certain function, and all the above transformations are invertible. We leave the details. If the
eigenvalue has multiplicity higher than 1, then these $\lambda-$forms span a
finite dimensional $\mathbb{C}-$vector space. Note that the proof of
Proposition \ref{prop1} implies $|\lambda|<0.342014_{+}$. Finally, the functional
equation ($\ref{eigen}$) gives the analytic continuation of $G_{\lambda}(z)$ to
the half-plane $\Re z\leq1$. Further, if $z\in\mathcal{U}$, where
$\mathcal{U}=\{0\leq \Re z\leq1\}\setminus\{|z|<1\}$, we can continue
$G_{\lambda}(z)$ to the region $\mathcal{U}+1$, and, inductively, to
$\mathcal{U}+n$, $n\in\mathbb{N}$. Let $\mathcal{U}_{0}$ be the union of these.
We can, obviously, continue $G_{\lambda}(z)$ to the set
$\mathcal{U}^{-1}_{0}+n$, $n\in\mathbb{N}$. Similar
iterations cover the described domain. $\blacksquare$ \\

\noindent Note that, in contrast to $G(z)$, we do not have a symmetry
property for $G_{\lambda}(z)$. \\

Next calculations produce the first few eigenvalues. Let the Taylor expansion
of $G_{\lambda}(z)$ is given by
\[
G_{\lambda}(z)=\sum\limits_{L=1}^{\infty}m_{L}^{(\lambda)}z^{L-1}.
\]
It converges in the unit circle, including its boundary (as is clear from
(\ref{eigen}), there exist all left derivatives at $z=1$). Thus,
$m_{L}^{(\lambda)}$ have the same vanishing properties as $m_{L}$ (which
guarantees the convergence of the series in ($\ref{rys}$)). And therefore, as
in the Proposition \ref{prop5}, we obtain:
\[
\lambda
m_{s}^{(\lambda)}=\sum\limits_{L=1}^{\infty}
(-1)^{L-1}c_{L+s}\binom{L+s-1}{s-1}m_{L}^{(\lambda)},\quad
s\geq 1.
\]
Here $c_{L}=\sum_{n=1}^{\infty}\frac{1}{2^{n}n^{L}}$. If we denote
$e_{s,L}=(-1)^{L-1}c_{L+s}\binom{L+s-1}{s-1}$, then $\lambda$ is the eigenvalue
of the infinite matrix $\mathcal{E}_{s,L=1}^{\infty}$. The numerical
calculations with the augmentation of this matrix at sufficiently high level
give the following first eigenvalues in decreasing order, with all digits
exact:
\[
\lambda_{1}=0.25553210_{+},\quad\lambda_{2}=-0.08892666_{+},\quad
\lambda_{3}=0.03261586_{+},\quad\lambda_{4}=-0.01217621_{+}.
\]
Paper \cite{alkauskas1} contains graphs of $G_{\lambda}(z)$ for the first six
eigenvalues, and, more importantly, ``pair-correlation'' results among
different eigenvalues, including the ``eigenvalue'' -1, which corresponds
precisely to $G(z)$. These results reveal the importance of $G_{\lambda}(z)$ in
the study of $F(x)$.
\section{$p-$adic distribution}

In the previous sections, we were interested in the distribution of the $n$th
generation of the tree $\mathcal{T}$ in the field of real numbers. Since the
set of non-equivalent valuations of $\mathbb{Q}$ contains a valuation
associated with any prime number $p$, it is natural to consider the
distribution of the set of each generation in the field of $p-$adic numbers
$\mathbb{Q}_{p}$. In this case we have an ultrametric inequality, which implies
that two circles are either co-centric or do not intersect. We define
\begin{eqnarray*}
F_{n}(z,\nu)=2^{-n+1}\#\{\frac{a}{b}\in\mathcal{T}^{(n)}:\text{ord}_{p}(\frac{a}{b}-z)\geq\nu\},\quad
z\in\mathbb{Q}_{p},\quad\nu\in\mathbb{Z}.
\end{eqnarray*}
(When $p$ is fixed, the subscript $p$ in $F_{n}$ is omitted). Note that in
order to calculate $F_{n}(z,\nu)$
we can confine to the case
$\text{ord}_{p}(z)<\nu$; otherwise
$\text{ord}_{p}(\frac{a}{b}-z)\geq\nu\Leftrightarrow\text{ord}_{p}(\frac{a}{b})\geq\nu$.
We shall calculate the limit distribution
$\mu_{p}(z,\nu)=\lim_{n\rightarrow\infty}F_{n}(z,\nu)$, and also some
characteristics of it, e.g. the zeta function
$$
Z_{p}(s)=\int\limits_{u\in\mathbb{Q}_{p}}|u|^{s}d\mu_{p},\quad s\in\mathbb{C},\quad z\in\mathbb{Q}_{p},
$$
where $|*|$ stands for the $p-$adic valuation.
\par

To illustrate how the method works, we will calculate the value of $F_{n}$ in two special cases. Let $p=2$ and
let $E(n)$ be the number of rational numbers in the $n$th generation with one of $a$ or $b$ being even, and let
$O(n)$ be the corresponding number fractions with both $a$ and $b$ odd. Then $E(n)+O(n)=2^{n-1}$. Since $\frac{a}{b}$
in the $n$th generation generates $\frac{a}{a+b}$ and $\frac{a+b}{b}$ in the $(n+1)$st generation, each fraction
$\frac{a}{b}$ with one of the $a$, $b$ even will generate one fraction with both numerator and denominator odd.
If both $a$, $b$ are odd, then their two offsprings will not be of this kind. Therefore, $O(n+1)=E(n)$. Similarly,
 $E(n+1)=E(n)+2O(n)$. This gives the recurrence $E(n+1)=E(n)+2E(n-1)$, $n\geq 2$, and this implies
$$
E(n)=\frac{2^{n}+2(-1)^{n}}{3},\quad O(n)=\frac{2^{n-1}+2(-1)^{n-1}}{3}, \quad \mu_{2}(0,0)=\frac{2}{3}.
$$
(For the last equality note that $\frac{a}{b}$ and $\frac{b}{a}$ simultaneously belong to $\mathcal{T}^{(n)}$,
and so the number of fractions with $\text{ord}_{2}(*)>0$ is $E(n)/2$). We will generalize this example to odd
prime $p\geq 3$. Let $L_{i}(n)$ be the part of the fractions in the $n$th generations such that $ab^{-1}\equiv i\bmod\,p$
for $0\leq i\leq p-1$ or $i=\infty$ (that is, $b\equiv0\bmod\,p$). Thus,
$$
\sum_{i\in\mathbb{F}_{p}\cup\infty}L_{i}(n)=1;
$$
in other words, $L_{i}(n)=F_{n}(i,1)$. For our later investigations we need a result from the theory of finite Markov chains.

\begin{lem}
Let $\textbf{A}$ be a matrix of a finite Markov chain with $s$ stages. That is,
$a_{i,j}\geq0$, and $\sum_{j=1}^{s}a_{i,j}=1$ for all $i$. Suppose that
$\mathbf{A}$ is irreducible (for all pairs $(i,j)$, and some $m$, the entry
$a_{i,j}^{(m)}$ of the matrix $\mathbf{A}^{m}$ is strictly positive), acyclic
and recurrent (this is satisfied, if all entries of $\mathbf{A}^{m}$ are
strictly positive for some $m$). Then the eigenvalue $1$ is simple and if
$\lambda$ is another eigenvalue, then $|\lambda|<1$, and $\mathbf{A}^{m}$, as
$m\rightarrow\infty$, tends to the matrix $\mathbf{B}$, with entries
$b_{i,j}=\pi_{j}$, where $(\pi_{1},...,\pi_{s})$ is a unique left eigenvector
with eigenvalue $1$, such that $\sum_{j=1}^{s}\pi_{j}=1$. \label{lem3}
\end{lem}

\noindent A proof of this lemma can be found in \cite{karl}, Section 3.1.,
Theorem 1.3.

\begin{prop}
$\mu_{p}(z,1)=\frac{1}{p+1}$ for $z\in\mathbb{Z}_{p}$. \label{prop6}
\end{prop}

{\it Proof. }Similarly as in the above example, a fraction $\frac{a}{b}$ from
the $n$th generation generates $\frac{a}{a+b}$ and $\frac{a+b}{b}$ in the
$(n+1)$st generation, and it is routine to check that
\begin{eqnarray}
L_{i}(n+1)=\frac{1}{2}L_{\frac{i}{1-i}}(n)+\frac{1}{2}L_{i-1}(n)\quad\mbox{for}\quad
i\in\mathbb{F}_{p}\cup\{\infty\}\label{padic},
\end{eqnarray}
(Here we make a natural convention for $\frac{i}{1-i}$ and $i-1$, if $i=1$ or $\infty$). In this equation,
it can happen that $i-1\equiv\frac{i}{1-i}\bmod\,p$; thus, $(2i-1)^{2}\equiv -3\bmod\,p$. The recurrence for
this particular $i$ is to be understood in the obvious way, $L_{i}(n+1)=L_{i-1}(n)$.  Therefore, if we denote
 the vector-column $(L_{\infty}(n),L_{0}(n),...,L_{p-1}(n))^{T}$ by $\mathbf{v}_{n}$, and if $\mathcal{A}$ is
 a matrix of the system $(\ref{padic})$, then $\mathbf{v}_{n+1}=\mathcal{A}\mathbf{v}_{n}$, and hence
$$
\mathbf{v}_{n}=\mathcal{A}^{n-1}\mathbf{v}_{1},
$$
where $\mathbf{v}_{1}=(0,0,1,0,...,0)^{T}$. In any particular case, this allows us two find the values
of $L_{i}$ explicitly. For example, if $p=7$, the characteristic polynomial is
$$
f(x)=\frac{1}{16}(x-1)(2x-1)(2x^{2}+1)(4x^{4}+2x^{3}+2x+1).
$$
The list of roots is
$$
\alpha_{1}=1,\quad\alpha=\frac{1}{2},\quad\alpha_{3,4}=\pm\frac{i}{\sqrt{2}},
\quad\alpha_{5,6,7,8}=\frac{-1-\sqrt{17}}{8}\pm\frac{\sqrt{1+\sqrt{17}}}{2\sqrt{2}},
$$
(with respect to the two values for the root $\sqrt{17}$), the matrix is
diagonalisible, and the Jordan normal form gives the expression
$$
L_{i}(n)=\sum_{s=1}^{8}C_{i,s}\alpha_{s}^{n}.
$$
Note that the elements in each row of the $(p+1)\times(p+1)$
matrix $\mathcal{A}$ are non-negative and sum up to $1$, and thus,
we have a matrix of a finite Markov chain. We need to check that
it is acyclic. Let $\tau(i)=i-1$, and $\sigma(i)=\frac{i}{1-i}$
for $i\in\mathbb{F}_{p}\cup\{\infty\}$. The entry $a_{i,j}^{(m)}$
of $\mathcal{A}^{m}$ is
$$
a_{i,j}^{(m)}=\sum_{i_{1},...,i_{m-1}}a_{i,i_{1}}\cdot a_{i_{1},i_{2}}\cdot...\cdot a_{i_{m-1},j}.
$$
Therefore, we need to check that for some fixed $m$, the
composition of $m$ $\sigma's$ or $\tau's$ leads from any $i$ to
any $j$. One checks directly that for any positive $k$, and
$i,j\in\mathbb{F}_{p}$,
\begin{eqnarray*}
\tau^{p-1-j}\circ\sigma\circ\tau^{k}\circ\sigma\circ\tau^{i-1}(i)&=&j,\\
\tau^{p-1-j}\circ\sigma\circ\tau^{k}(\infty)&=&j,\\
\tau^{k}\circ\sigma\circ\tau^{i-1}(i)&=&\infty;
\end{eqnarray*}
(for $i=0$, we write $\tau^{-1}$ for $\tau^{p-1}$). For each pair $(i,j)$, choose $k$ in
order the amount of compositions used to be equal (say, to $m$). Then obviously all entries
 of $\mathcal{A}^{m}$ are positive, ant this matrix satisfies the conditions of Lemma \ref{lem3}. Since
 all columns also sum up to $1$, $(\pi_{1},...,\pi_{p+1})$, $\pi_{j}=\frac{1}{p+1}$, $1\leq j\leq p+1$,
 is the needed eigenvector. This proves the Proposition. $\blacksquare$\\

Next Theorem describes $\mu(z,\nu)$ in all cases.
\begin{thm}
Let $\nu\in\mathbb{Z}$ and $z\in\mathbb{Q}_{p}$, and $\text{ord}_{p}(z)<\nu$ (or $z=0$). Then, if $z$ is $p-$adic integer,
$$
\mu(z,\nu)=\frac{1}{p^{\nu}+p^{\nu-1}}.
$$
If $z$ is not integer, $\text{ord}_{p}(z)=-\lambda<0$,
$$
\mu(z,\nu)=\frac{1}{p^{\nu+2\lambda}+p^{\nu+2\lambda-1}}.
$$
For $z=0$, $-\nu\leq0$, we have
$$
\mu(0,-\nu)=1-\frac{1}{p^{\nu+1}+p^{\nu}}.
$$
\label{thm4}
\end{thm}

\noindent This Theorem allows the computation of the associated zeta-function:

\begin{cor} For $s$ in the strip $-1<\Re{s}<1$,
\begin{eqnarray*}
Z_{p}(s)=\int_{u\in\mathbb{Q}_{p}}|u|^{s}d\mu_{p}=
\frac{(p-1)^{2}}{(p-p^{-s})(p-p^{s})},
\end{eqnarray*}
and $Z_{p}(s)=Z_{p}(-s)$.
\end{cor}

\noindent The proof is straightforward. It should be noted that this expression encodes all the values of $\mu(0,\nu)$ for $\nu\in\mathbb{Z}$.
\par\bigskip

{\it Proof of Theorem \ref{thm4}. }For shortness, when $p$ is fixed, denote
$\text{ord}_{p}(*)$ by $v(*)$. As before, we want a recurrence relation among
the numbers $F_{n}(i,\kappa)$, $i\in\mathbb{Q}_{+}$. For each integral
$\kappa$, we can confine to the case $i<p^{\kappa}$. If $i=0$, we only consider
$\kappa>0$ and call these pairs $(i,\kappa)$ ``admissible''. We also include
$G_{n}(0,-\kappa)$ for $\kappa\geq 1$, where these values are defined in the
same manner as $F_{n}$, only inverting the inequality, considering
$\frac{a}{b}\in\mathcal{T}^{(n)}$, such that $v(\frac{a}{b})\leq -\kappa$; this
is the ratio of fractions in the $n$th generation outside this circle. As
before, a fraction $\frac{a}{b}$ in the $n$th generation generates the
fractions $\frac{a}{a+b}$ and $\frac{a+b}{b}$ in the $(n+1)$st generation. Let
$\tau(i,\kappa)=((i-1)\bmod\,p^{\kappa},\kappa)$. Then for all admissible pairs
$(i,\kappa)$, $i\neq0$, the pair $\tau(i,\kappa)$ is also admissible, and
$$
v(\frac{a+b}{b}-i)=\kappa\Leftrightarrow v(\frac{a}{b}-(i-1))=\kappa.
$$
Second, if $\frac{a}{a+b}=i+p^{\kappa}u$, $i\neq 1$,
$u\in\mathbb{Z}_{p}$, and $(i,\kappa)$ is admissible, then
$$
\frac{a}{b}-\frac{i}{1-i}=\frac{p^{\kappa}u}{(1-i)(1-i-p^{\kappa}u)}.
$$
Since $v(\frac{i}{1-i})=v(i)-v(1-i)$, this is $0$ unless $i$ is an
integer, equals to $v(i)$ if the latter is $>0$ and equals to
$-v(1-i)$ if $v(1-i)>0$. Further, this difference has valuation
$\geq\kappa_{0}=\kappa$, if $i\in\mathbb{Z},i\not\equiv
1\bmod\,p$, valuation $\geq\kappa_{0}=\kappa-2v(1-i)$, if
$i\in\mathbb{Z},i\equiv1\bmod\,p$, and valuation
$\geq\kappa_{0}=\kappa-2v(i)$ if $i$ is not integer. In all three
cases, easy to check, that, if we define
$i_{0}=\frac{i}{1-i}\bmod\,p^{\kappa_{0}}$, the pair
$\sigma(i,\kappa)\mathop{=}^{\text{def}}(i_{0},\kappa_{0})$ is admissible.
For the converse, let $\frac{a}{b}=i_{0}+p^{\kappa_{0}}u$,
$u\in\mathbb{Z}_{p}$. Then
$$
\frac{a}{a+b}-\frac{i_{0}}{1+i_{0}}=\frac{p^{\kappa_{0}}}{(1+i_{0}+p^{\kappa_{0}}u)(1+i_{0})}.
$$
If $i=\frac{i_{0}}{1+i_{0}}$ is a $p-$adic integer, $i\not\equiv 1\bmod\,p$, this has a valuation $\geq\kappa=\kappa_{0}$; if $i$ is a $p-$adic
integer, $i\equiv 1(p)$, this has valuation
$
\geq\kappa=\kappa_{0}-2v(i_{0})=\kappa_{0}+2v(1-i);
$
if $i$ is not a $p-$adic integer, this has valuation
$
\geq\kappa=\kappa_{0}-2v(1+i_{0})=\kappa_{0}+2v(i).
$
Thus,
$$
v(\frac{a}{a+b}-i)\geq\kappa\Leftrightarrow v(\frac{a}{b}-i_{0})\geq\kappa_{0}.
$$
Let $i=1$. If $\frac{a}{a+b}=1+p^{\kappa}u$, then $\kappa>0$,
$u\in\mathbb{Z}_{p}$, and we obtain
$\frac{a}{b}=-1-\frac{1}{p^{\kappa}u}$,
$v(\frac{a}{b})\leq-\kappa$. Converse is also true. Finally, for
$\kappa\geq1$,
$$
v(\frac{a+b}{b})\leq-\kappa\Leftrightarrow v(\frac{a}{b})\leq-\kappa,
$$
and
$$
v(\frac{a}{a+b})\leq-\kappa\Leftrightarrow v(\frac{a}{b}+1)\geq\kappa.
$$
Therefore, we have the recurrence relations:
\begin{eqnarray}\left\{ \begin{array}{l@{\,}l}
F_{n+1}(i,\kappa)&=\frac{1}{2}F_{n}(\tau(i,\kappa))+\frac{1}{2}F_{n}(\sigma(i,\kappa)),
\text{ if }
(i,\kappa)\text{ is admissible},\\
F_{n+1}(1,\kappa)&=\frac{1}{2}F_{n}(0,\kappa)+\frac{1}{2}G_{n}(0,-\kappa), \kappa\geq1,\\
G_{n+1}(0,-\kappa)&=\frac{1}{2}G_{n}(0,-\kappa)+\frac{1}{2}F_{n}(-1,\kappa), \kappa\geq1.\\
\end{array} \right.\label{pad}
\end{eqnarray}
Thus, we have an infinite matrix $\mathcal{A}$, which is a change
matrix for the Markov chain. If $\mathbf{v}_{n}$ is an infinite
vector-column of $F_{n}'$s and $G_{n}'$s, then
$\mathbf{v}_{n+1}=\mathcal{A}\mathbf{v}_{n}$, and, as before,
$\mathbf{v}_{n}=\mathcal{A}^{n-1}\mathbf{v}_{1}$. It is direct to
check that each column also contains exactly two nonzero entries
$\frac{1}{2}$, or one entry, equal to $1$. In terms of Markov
chains, we need to determine the classes of orbits. Then in proper
rearranging, the matrix $\mathcal{A}$ looks like
\begin{eqnarray*}
\begin{pmatrix}
\mathbf{P}_{1} & 0 & \dots & 0 & \dots\\
0              &\mathbf{P}_{2} & \dots & 0 & \dots\\
\vdots &  & \ddots & \vdots& \vdots\\
0 & 0 & \dots & \mathbf{P}_{s}& 0\\
\vdots & \vdots & \dots & 0 & \ddots
\end{pmatrix},
\end{eqnarray*}
where $\mathbf{P}_{s}$ are finite Markov matrices. Thus, we claim
that the length of each orbit is finite, every orbit has a
representative $G_{*}(0,-\kappa)$, $\kappa\geq 1$, the length of
it is $p^{\kappa}+p^{\kappa-1}$, and the matrix is recurrent (that
is, every two positions communicate). In fact, from the system
above and form the expression of the maps $\tau(i,\kappa)$ and
$\sigma(i,\kappa)$, the direct check shows that the complete list
of the orbit of $G_{*}(0,-\kappa)$ consists of (and each pair of
states are communicating):
\begin{eqnarray*}
&G_{*}(0,-\kappa),\\
&F_{*}(i,\kappa)\quad (i=0,1,2,...,p^{\kappa}-1),\\
&F_{*}(p^{-\lambda}u,\kappa-2\lambda)\quad (\lambda=1,2,...,\kappa-1, u\in\mathbb{N},u\not\equiv 0 \bmod\,p,
u\leq p^{\kappa-\lambda}).
\end{eqnarray*}
In total, we have
$$
1+p^{\kappa}+\sum\limits_{\lambda=1}^{\kappa-1}(p^{\kappa-\lambda}-p^{\kappa-\lambda-1})= p^{\kappa}+p^{\kappa-1}
$$
members in the orbit. Thus, each $\mathbf{P}_{\kappa}$ in the matrix above is a
finite dimensional $\ell_{\kappa}\times\ell_{\kappa}$ matrix, where
$\ell_{\kappa}=p^{\kappa}+p^{\kappa-1}$. For $\kappa=1$, the matrix
$\mathbf{P}_{1}$ is exactly the matrix of the system (\ref{padic}). As noted
above, the vector column $(1,1,...,1)^{T}$ is the left eigenvector. As in the
previous Proposition, it is straightforward to check that this matrix is
irreducible and acyclic (that is, the entries of $\mathbf{P}_{\kappa}^{n}$ are
strictly positive for sufficiently large $n$). In fact, since by our
observation, each two members in the orbit communicate, and since we have a
move $G_{*}(0,-\kappa)\rightarrow G_{*}(0,-\kappa)$, the proof of the last
statement is immediate: there exists $n$ such that any position is reachable
from another in exactly $n$ moves, and this can be achieved at the expense of
the move just described. Therefore, all entries of $\mathbf{P}_{\kappa}^{n}$
are strictly positive. Thus, the claim of the Theorem follows from the Lemma
\ref{lem3}. $\blacksquare$

\section{Conclusion}
We end the paper with the following remarks. As is implied by Theorem
\ref{thm4}, the measure $\mu_{p}$ of those rationals in the Calkin-Wilf tree
which are invertible elements of $\mathbb{Z}_{p}$ is equal to
$\frac{p-1}{p+1}$. We follow the line of the Tate thesis \cite{cassels}, and modify
this measure in order $\mathbb{Z}_{p}^{*}$ to have measure $1$; accordingly,
let us define $\mu'_{p}=\frac{p+1}{p-1}\mu_{p}$. Thus, we are lead to the
formal definition of the zeta function
\[
\zeta_{\mathcal{T}}(s)=\prod_{p}\int_{u\in\mathbb{Q}_{p}}|u|^{s}\d\mu_{p}'=
\prod_{p}\Big{(}1-\frac{1}{p^{2}}\Big{)}
\prod_{p}\frac{1}{1-p^{-s-1}}\cdot\frac{1}{1-p^{s-1}}=
\frac{6}{\pi^{2}}\zeta(s+1)\zeta(-s+1).
\]
This product diverges everywhere; nevertheless, if we apply the
functional equation of the Riemann $\zeta$ function for the second
multiplier, we obtain
\[
\zeta_{\mathcal{T}}(s)=\frac{12}{\pi^{2}}(2\pi)^{-s}\cos\Big{(}\frac{\pi
s}{2}\Big{)}\Gamma(s)\zeta(s)\zeta(s+1).
\]
 From the above definition it is clear that, formally, this zeta function is the sum of the form
$\sum_{r\in\mathbb{Q}^{+}}\mu_{r}r^{-s}$, where, if $r\in\mathbb{Q}_{+}$, and $\mu_{r}$ stands for the limit measure of those rationals in the $n$th
generation of $\mathcal{T}$, which have precisely the same valuation as $r$ at
every prime which appears in the decomposition of $r$, times the factor
$\prod_{\text{ ord}_{p}(r)\neq0}\frac{p+1}{p-1}$. Surprisingly, the product
$\zeta(s)\zeta(s+1)$ is the zeta function of the Eisenstein series $G_{1}(z)$,
which is related to the distribution of rationals in $\mathcal{T}$ at the
infinite prime $\mathbb{Q}_{\infty}=\mathbb{R}$. In fact,
\[
\int\limits_{0}^{\infty}\Big{(}G_{1}(iz)-G_{1}(i\infty)\Big{)}z^{s-1}\d z=
-8\pi^{2}(2\pi)^{-s}\Gamma(s)\zeta(s)\zeta(s+1).
\]
This is a strong motivation to investigate the tree $\mathcal{T}$ and the
Minkowski question mark function in a more general - idelic - setting, thus
revealing the true connection between $p-$adic and real distribution, and
clarifying the nature of continued fractions in this direction. We hope to
implement this in the subsequent papers.\\
\indent Unfortunately, currently we left the most interesting question (the explicit
description of the moments of $?(x)$) unanswered. It is desirable to
give the function $G(z)$ and dyadic forms $G_{\lambda}(z)$ certain other
description than the one which arises directly from the tree $\mathcal{T}$.
This is in part accomplished in \cite{alkauskas1} and \cite{alkauskas3}. These
papers form a direct continuation of the current work. Among the other results, the dyadic zeta function $\zeta_{\mathcal{M}}(s)$ is introduced: it is given by
$\zeta_{\mathcal{M}}(s)\Gamma(s+1)=\int_{0}^{\infty}x^{s}\d F(x)$; the nature
of dyadic eigenfunctions $G_{\lambda}(z)$ is clarified; certain integrals
which involve $F(x)$ are computed; and finally, this research culminates with the proof that in the half plane $\Re z<1$ the dyadic period function $G(z)$ can be represented as an
absolutely convergent series of rational functions with rational coefficients.
Possibly, this technique can find its applications in the study of period
function for Maass wave forms.

\par\bigskip

\noindent Giedrius Alkauskas, School of Mathematical Sciences, University of
Nottingham, University Park, Nottingham NG7 2RD
United Kingdom\\
 {\tt giedrius.alkauskas@maths.nottingham.ac.uk}
\end{document}